\documentclass[10pt]{article}

\usepackage[hidelinks]{hyperref}
\usepackage{amsmath}
\usepackage{eqnarray,amsmath}
\usepackage{amssymb}
\usepackage{amsfonts}
\usepackage[usenames,dvipsnames]{pstricks}
\usepackage{graphicx}
\usepackage{subcaption}
\usepackage[labelformat=parens,labelsep=quad,skip=3pt]{caption}
\usepackage{qtree}
\usepackage{algorithm}
\usepackage[noend]{algpseudocode}
\usepackage[utf8]{inputenc}
\usepackage{amsthm}
\usepackage{pstricks-add,}
\usepackage{mathtools}
\usepackage{dirtytalk}

\usepackage{xcolor}
\usepackage{hyperref}

\usepackage[english]{babel} % needed for "blindtext"
\usepackage[pangram]{blindtext}

\usepackage{array}
\usepackage{rotating}
\usepackage{tikz}
\usepackage{blkarray}

\usepackage{MnSymbol} % It is for \udots command

\usepackage{bbm}

\usepackage{nicematrix}
\NiceMatrixOptions{
code-for-first-row = \color{blue} ,
code-for-last-row = \color{blue} ,
code-for-first-col = \color{blue} ,
code-for-last-col = \color{blue}
}

\hypersetup{colorlinks=false,linkbordercolor=red,linkcolor=green,pdfborderstyle={/S/U/W 1}}

\newtheorem{theorem}{Theorem}[section]

\newtheorem{remark}[theorem]{Remark}

\newtheorem{example}[theorem]{Example}
\newtheorem{definition}[theorem]{Definition}
\tikzset{use/.style={}}
\tikzset{rectangle/.append style={circle, fill=blue!25}}

\newcommand{\PosetInsertion}[3]{
    \ensuremath{
    \mathop{
    \begin{tikzpicture}[scale=0.25]
        \ifnum#1=1 \fill[black!75] (0,0.5) rectangle (0.5,1); \fi
        \ifnum#2=1 \fill[black!75] (0,0) rectangle (0.5,0.5); \fi
        \ifnum#3=1 \fill[black!75] (0.5,0) rectangle (1,0.5); \fi
        \draw[thin](0,0.5) rectangle (0.5,1);
        \draw[thin](0,0) rectangle (0.5,0.5);
        \draw[thin](0.5,0) rectangle (1,0.5);
    \end{tikzpicture}
    }}
}

\newcommand{\PosetComposition}[2]{
    \,
    \begin{tikzpicture}[scale=0.25]
        \ifnum#1=1 \fill[black!75] (0,0.5) rectangle (0.5,1); \fi
        \ifnum#2=1 \fill[black!75] (0.5,0) rectangle (1,0.5); \fi
        \draw[thin](0,0.5) rectangle (0.5,1);
        \draw[thin](0.5,0) rectangle (1,0.5);
    \end{tikzpicture}
    \,
}

\begin{document}

\title{Poset Matrix Structure Via Partial Composition Operations}
\date{}
\author{Arnauld Mesinga Mwafise\\
{\footnotesize $^{}$ \textit{}}\\
{\footnotesize mesingamwafisea@gmail.com}} \maketitle

\providecommand{\msc}[1]
{
  \small	
  \textbf{\textit{MSC classes:---}} #1

}

\providecommand{\keywords}[1]
{
  \small	
  \textbf{\textit{Keywords---}} #1

}

\begin{abstract}

This paper examines the structure of poset matrices by formulating a set of new construction rules for this purpose. In this direction, the technique of partial composition operation will be introduced as the basis for the construction of poset matrices of any given size by extending the combinatorial setting of species of structures to poset matrices. More specifically, three new partial composition operations that apply to poset matrices are defined as the foundation for this study. Several new structural properties derived from viewing any poset matrix and its dual in terms of these operations are highlighted.

\end{abstract}

\msc{06A06, 06A07, 15A36}
\keywords{poset matrices, partial composition operations}

\section{Introduction}

A nonempty set $X$ of size $n$ which is endowed with the order relation $\leqslant$ is called a partially ordered set or simply a poset,
 denoted by ${\mathcal P}=(X,\leqslant)$, , if the relation $\leqslant$ is 
 \begin{enumerate}
  \item reflexive: $x\in X, x\leq x$, 
  \item  antisymmetric: $x,y\in X$ if $x\leq y$ and $y\leq x\implies x=y$ 
  \item transitive: $x,y,z\in X$ if $x\leq y$ and $y\leq z\implies x\leq z.$
 \end{enumerate}
 In general, all posets can be structurally represented using directed acyclic graphs in which there is a path from $x$ to $y$ whenever $x\leq y$ and with a transitively reduced graph connection path to $x\leq z$  whenever $x\leq y$ and $y\leq z.$
For any labeled n-poset,  there exist a corresponding matrix representation that encodes the partial order relationships between the elements in the set of labels. This matrix is called the labeled matrix of a labeled n- poset or simply a \textit{poset matrix}. Therefore, every finite poset $\mathcal{P}$ can be associated  with an $n\times n$ binary matrix $A=[a_{i,j}]$ labeled by elements of the set $X$ if the following three conditions are satisfied for all $i,j,k\in X$:
 
\begin{enumerate}
\item[{\rm (i)}] $A$ is reflexive, i.e. $a_{ii}=1;$ 
\item[{\rm (ii)}] $A$ is antisymmetric, i.e. if $a_{ij}=1$ then $a_{ji}=0$;
\item[{\rm (iii)}] $A$ is transitive, i.e. if $a_{ij}=1$ and $a_{jk}=1$ then $a_{ik}=1$.
\end{enumerate}

\noindent A fundamental understanding about the symmetries and structural arrangement of the $(0,1)$ entries of any  poset matrix can be highlighted by viewing it in terms of smaller connected or disconnected subposet matrices. The outcome of this process leads to an intuition on the enumeration of poset matrices. Several operations have been introduced to connect two or more poset matrices to create larger poset matrices as can be seen with the specific case of series parallel poset matrices using the ordinal sum and the direct sum operations \cite{mohammad, riordanp}. Furthermore, the imporatnace of developing new techniques for constructing new posets from old ones has been highlighted  as a key research objective in poset theory \cite{stanley}.   The partial composition operations in the next section provides a more generalized technique for constructing any poset matrix of a given size from smaller ones. The origins of partial composition operations can be traced back to the study of the combinatorics of labeled rooted trees, graphs and words in the context of operads \cite{aval,chapoton,joyal,bergeron}. More recently, it has been applied to posets represented by their Hasse diagrams \cite{fauvet}. In particular, the substitution operations performed by the grafting or insertion at specified labeled points of these combinatorial structures forms the basis of the underlying partial composition operations based on specific rules that apply to each combinatorial structure. Moreover, the most important consideration to defining the operad structure of a combinatorial object is  to determine an appropriate partial composition operation which specifies a set of rules on how to connect two or more  labeled combinatorial structures to form a new combinatorial structure of that same type satisfying an identity, parallel and nested associativity axioms \cite{chapoton}. Operads are also considered a special kind of species. Combinatorial species was introduced by Andr$\acute{\text{e}}$ Joyal \cite{joyal} as  a novel way to analyze complicated structures  by describing them in terms of transformations and combinations of simpler structures. A key aspect of species is the way it associates structures to some base set of labels from which many possible structures can be formed. Generally, it provides a universal way of dealing with generating series in a monoidal categorical setting  \cite{mendez}.
 Similar to labeled posets represented by Hasse diagrams of their directed graphs which  can be viewed in terms of combinatorial species,  this can also be extended to their corresponding labeled poset matrices as presented in this study. 

The rest of the paper is structured as follows. In section two, the three partial composition operation rules for poset matrix construction are formulated. In section three, we provide some key structural properties of poset matrices and their dual representations based on partial composition operations.

\section{Partial composition operations of poset matrices}

Recall that for a poset ${\mathcal P}=(X,\leqslant)$ the \textit{minimal element} refers to any element $p\in (X,\leqslant)$ if there is no other element $q\in (X,\leqslant)$ such that $p<q$ . In the Hasse diagram, it can be identified as an element with no incoming(downward) edges. On the other hand, a \textit{maximal element} is an element $p\in (X,\leqslant)$ if there  is no other element $q\in (X,\leqslant)$ such that $q<p.$    It can be identified in the Hasse diagram of the poset as the element with no outgoing(forward) edge. We provide analogous results for identifying the  maximal and minimal elements of labeled poset matrices in the theorem that follows. 

\begin{theorem}\label{minmax}
Let $A=[a_{i,j}]$ be a poset matrix of size $n$ and let $L$ be any $n$-element labeling set for the row and column indices of $A.$ Then whenever $i\neq j$,
\begin{enumerate} 
\item the set of all maximal elements of $A$ consist of all $j\in L$ such that $a_{i,j}=0$ for every $i\in L.$
\item the set of all minimal elements of $A$ consist of all $i\in L$ such that $a_{i,j}=0$ for every $j\in L.$
\end{enumerate} 

\end{theorem}

\begin{remark}
For a poset matrix $A$ we denote by $max(A)$  the set of all maximal elements of $A$ and  $min(A)$  the set of all minimal elements of $A.$  
\end{remark}

\begin{example}\label{matuse}
\end{example}

\noindent Consider the poset matrices:

$\begin{pNiceMatrix}[first-row,first-col]
&5 & 6&7\\
5&1 & 0&0 \\
6&1 & 1&0 \\
7&1 & 0&1  \\
\end{pNiceMatrix}=B\qquad \begin{pNiceMatrix}[first-row,first-col]
&1&2 & 3 &4\\
 1&1 & 0 & 0&0 \\
 2&1 & 1 & 0&0 \\
 3&1 & 1 & 1&0 \\
 4&1 & 1 & 0&1 \\
\end{pNiceMatrix}=A.$

Using Theorem \ref{minmax},  the label index of the minimal element of  $A:=\text{min}(A)=\{1\}$ and the labeled indices of the maximal element of $A:=\text{max}(A)=\{3,4\}.$ On the other hand, the label index of the minimal element of  $B:=\text{min}(B)=\{5\}$ and the label indices of the maximal element of $B:=\text{max}(B)=\{6,7\}.$

\noindent Consider two input poset matrices $A_n$ and $B_m$ of size $n$ and $m$ respectively. These two input poset matrices are concatenated using the procedure known as the \textit{partial composition operation} which is most often denoted by $\circ_i$ such that  $A_n\circ_i B_m=C_{n+m-1}$ where $i\in\{1,...,n\}$ is considered the row and column index removed from $A_n$ and replaced by indices(labels) from $B_m$ in this study. The output poset matrix $C_{n+m-1}$ of size $n+m-1$ resulting from the construction will have its entries determined by the three types of proposed partial composition operations in this study. For the purpose of clarity, we introduce the following matrix which depicts the output matrix from the partial composition operation.
\begin{equation}\label{mainmat}
C_{n+m-1}=\begin{pNiceMatrix}[first-row,first-col]
    & 1 & \dots  & i-1 & n+1 &\dots & n+m & i+1 &\dots & n \\
1 & a_{1,1}   & \dots   & a_{1,i-1}   & 0& 0 & 0   & 0   & 0& 0    \\
\vdots & \vdots   & \ddots  & \vdots   & \vdots &  \vdots & \vdots   & \vdots  & \vdots& \vdots \\
 i-1 & a_{i-1,1}   & \dots   & a_{i-1,i-1}   & 0 & 0  & 0   & 0& 0&0  \\
n+1 & u_{n+1,1}   & \dots   & u_{n+1,i-1}   & b_{n+1,n+1} & \dots   & b_{n+1,n+m}   & 0 & 0& 0  \\
\vdots & \vdots   & \ddots   & \vdots   & \vdots& \ddots & \vdots   & \vdots   &  \vdots & \vdots   \\
n+m& u_{n+m,1}   & \dots  & u_{n+m,i-1}   & b_{n+m,n+1}& \dots & b_{n+m,n+m}   & 0   & 0& 0  \\
i+1 & a_{i+1,1}   & \dots   & a_{i+1,i-1}   & v_{i+1,n+1} & \dots  & v_{i+1,n+m}   & a_{i+1,i+1}   & \dots & a_{i+1,n}  \\
\vdots & \vdots   & \ddots   & \vdots    & \vdots & \ddots  & \vdots   & \vdots   &\ddots& \vdots \\
n & a_{n,1}   & \dots   & a_{n,i-1}   & v_{n,n+1} & \dots  & v_{n,n+m}   & a_{n,i+1}   & \dots & a_{n,n}  \\
\end{pNiceMatrix}
\end{equation}

The matrix $C_{n+m-1}$ can also be viewed as follows:
\begin{eqnarray}\label{partial}
C_{n+m-1}=\begin{pmatrix}
A_{n_1}&\vline &O&\vline&O\\
\hline
U&\vline &B_m&\vline&O\\
\hline
A_{n_3}&\vline &V&\vline&A_{n_2}\\
\end{pmatrix}\;\;(1<i<n),
\end{eqnarray}
and
\begin{align*}
A_n\circ_1 B_m =
\left[\begin{array}{c|c}
B_m & O\\
\hline
V & A_{n_2}
\end{array}\right]\;\;{\rm and}\;\;
A_n\circ_n B_m =
\left[\begin{array}{c|c}
A_{n_1} & O \\
\hline
U & B_m
\end{array}\right],
\end{align*}
where $U$ and $V$ are (0,1)-matrices of sizes $m\times (i-1)$ and $(n-i)\times m$, respectively.

\noindent The $(0,1)$ entries of the submatrices $U$ and $V$ in the matrix $C_{n+m-1}$ can be determined with the purpose of obtaining a poset matrix using any of the following three methods as described below.
\\
\\
\noindent Let $A_n=(a_{y,z}),$ $U=(u_{y,z})$ and $V=(v_{y,z}).$ Then for $A_n\circ_i B_m=C_{n+m-1}$  we have that:

\begin{enumerate}

\item [{\rm(1)}]  The first partial composition operation  denoted  by $\scalebox{0.7}{$\square$}_i$ can be described as follows. 
 If $y$ and $z$ are the row and column labels of $U$ then $u_{y,z}=a_{i,z}$ . Otherwise, if $y$ and $z$ are the row and column labels of $V$ then $v_{y,z}=a_{y,i}.$  

\item[{\rm(2)}]  The second partial composition operation  denoted  by $\scalebox{0.7}{$\blacktriangledown$}_i$ can be described as follows.

Let $P$ denote the set of all minimal elements of the poset matrix  $B_m$. Then 

\begin{equation}\label{min}
	v_{y,z}=
	\begin{cases}
		a_{y,i} & \text{if z}\in P\\
		0 & \text{otherwise.}
	\end{cases}
\qquad\; \text{and}\;\;u_{y,z}=a_{i,z}.\end{equation}

\item[{\rm(3)}]  The third partial composition operation  denoted  by $\scalebox{0.7}{$\blacktriangle$}_i$  can be described as follows.

Let $R$ denote the set of all maximal elements of the poset matrix $B_m$. Then 

\begin{equation}\label{min}
	u_{y,z}=
	\begin{cases}
		a_{i,z} & \text{if  y}\in R \\
		0 & \text{otherwise.}
	\end{cases}
\qquad\; \text{and}\;\;v_{y,z}=a_{y,i}.\end{equation}
  
\end{enumerate}

\begin{example}
\end{example}

Consider the poset matrix

$$A=B=\begin{blockarray}{ccccc}
1&2 & 3 &4\\
\begin{block}{(cccc)c}
 1 & 0 & 0&0 &1\\
 1 & 1 & 0&0 &2\\
 1 & 0 & 1&0 &3\\
 1 & 1 & 1&1 &4\\
\end{block}
\end{blockarray}$$
By Theorem \ref{minmax},  $\text{min}(B)=\{1\}$ and $\text{max}(B)=\{4\}$. Then, we have the following output poset matrices:

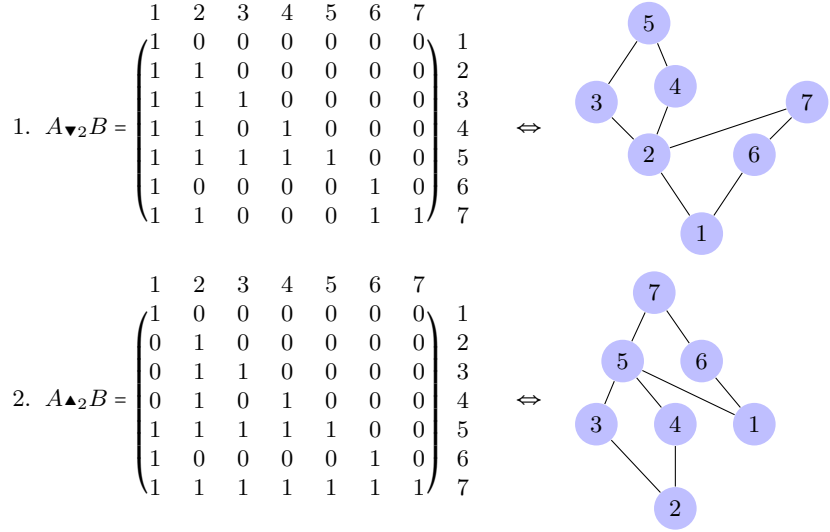
\begin{enumerate}

\item $A\scalebox{0.7}{$\blacktriangledown$}_2 B=
\begin{blockarray}{cccccccc}
1&2 & 3 &4&5&6&7\\
\begin{block}{(ccccccc)c}
  1 & 0 & 0 & 0 & 0 & 0 & 0&1 \\
1 & 1 & 0 & 0 & 0 & 0 & 0&2 \\
 1 & 1 & 1 & 0 & 0 & 0 & 0&3 \\
 1 & 1 & 0 & 1 & 0 & 0 & 0&4 \\
 1 & 1 & 1 & 1 & 1 & 0 & 0&5 \\
 1 & 0 & 0 & 0 & 0 & 1 & 0&6 \\
 1 & 1 & 0 & 0 & 0 & 1 & 1&7 \\
\end{block}
\end{blockarray}\quad\Leftrightarrow\quad$
\begin{minipage}[c]{.20\textwidth}
\begin{tikzpicture}
  [scale=.7,auto=center,every node/.style={circle,fill=blue!20}]
 \node[rectangle] (c1) at (9,1.5) {3};
  \node[rectangle] (c2) at (12,0.5) {6};
  \node[rectangle] (c4) at (10.5,1.8) {4};
   \node[rectangle] (c5) at (13,1.5) {7};
     \node[rectangle] (c6) at (10,0.5) {2};
       \node[rectangle] (c7) at (11,-1) {1};
        \node[rectangle] (c8) at (10,3) {5};
      \draw (c5) -- (c2);
         \draw (c2) -- (c7);
            \draw (c7) -- (c6);
      \draw (c1) -- (c6);
        \draw (c6) -- (c4);
         \draw (c6) -- (c5);
         \draw (c8) -- (c4);
         \draw (c8) -- (c1);

 \end{tikzpicture}
 
\end{minipage}

\item $A\scalebox{0.7}{$\blacktriangle$}_2 B=
\begin{blockarray}{cccccccc}
1&2 & 3 &4&5&6&7\\
\begin{block}{(ccccccc)c}
  1 & 0 & 0 & 0 & 0 & 0 & 0&1 \\
 0 & 1 & 0 & 0 & 0 & 0 & 0&2 \\
 0 & 1 & 1 & 0 & 0 & 0 & 0&3 \\
 0 & 1 & 0 & 1 & 0 & 0 & 0&4 \\
 1 & 1 & 1 & 1 & 1 & 0 & 0&5 \\
 1 & 0 & 0 & 0 & 0 & 1 & 0&6 \\
 1 & 1 & 1 & 1 & 1 & 1 & 1&7 \\
\end{block}
\end{blockarray}\quad\Leftrightarrow\quad$
\begin{minipage}[c]{.20\textwidth}
\begin{tikzpicture}
  [scale=.7,auto=center,every node/.style={circle,fill=blue!20}]
 \node[rectangle] (c1) at (10,1.7) {5};
  \node[rectangle] (c2) at (11,0.5) {4};
  \node[rectangle] (c4) at (11.5,1.7) {6};
   \node[rectangle] (c5) at (12.5,0.5) {1};
     \node[rectangle] (c6) at (9.5,0.5) {3};
       \node[rectangle] (c7) at (11,-1.1) {2};
        \node[rectangle] (c8) at (10.6,3) {7};
      \draw (c1) -- (c2);
         \draw (c1) -- (c5);
            \draw (c7) -- (c6);
      \draw (c1) -- (c6);
         \draw (c8) -- (c4);
         \draw (c8) -- (c1);
             \draw (c4) -- (c5);
                 \draw (c2) -- (c7);

 \end{tikzpicture}
 
\end{minipage}

\end{enumerate}

\begin{theorem}\label{subposets}
Given an $n$-element set $X=\{1,...,n\}$ of labels of the poset matrix $A_n$ of size $n,$ then every $m\times m$ submatrix of $A_n$ whose rows and columns are labeled from the same m-element subset of $X$ is a poset matrix.
\end{theorem}

\begin{theorem}\label{posetpartial}
Let $A_n$ and $B_m$ denote $n\times n$ and $m\times m$ poset matrices respectively.The matrix $C$ of size $(n+m-1)\times(n+m-1)$ obtained from $A_n\circ_i B_m,$ where $\circ_i\in\{\scalebox{0.7}{$\square$}_i,\scalebox{0.7}{$\blacktriangledown$}_i,\scalebox{0.7}{$\blacktriangle$}_i\}$ is a poset matrix for all $n,m\in\mathbb{N}.$
\end{theorem}

\section{Some key structural properties of poset matrices using their partial composition operations}

\noindent It is well known that a poset is connected if the Hasse diagram of its directed graph is connected. For the purpose of this work, we make the following definition categorizing poset matrices in terms of their connectivity with respect to their corresponding subposet matrix structures as follows.

\begin{definition}
Let $A=[a_{i,j}]$ be an $n\times n$ poset matrix. If there exist at least one subposet matrix of $A$ defined on a  non-empty subset $\alpha$ of $[n]$ such that for each $k\in \alpha,$  $a_{i,k}=0$ and $a_{k,j}=0$ whenever $i,j\notin\alpha,$ then we refer to $A$ as a \textit{disconnected poset matrix.} Otherwise, we say that $A$ is a \textit{connected poset matrix}.
\end{definition}

The following theorem characterizes disconnected poset matrices in terms of partial composition operation.

\begin{theorem}\label{dpm}
Let $\text{A}_n$ and $\text{B}_m$ be poset matrices of order $n$ and $m$ respectively. The partial composition operation $A_n\;\scalebox{0.7}{$\square$}_i\;B_m=C_{n+m-1}$ forms a disconnected poset matrix for all $i$ if and only if $A_n$ is a disconnected poset matrix.
\end{theorem}

In poset theory, for a poset $(P,\leq)$ if $x,y\in P$ such that $x\leq y$ then the dual poset matrix $(P^{*},\leq)$ with $x,y\in P^{*}$ satisfies the reverse order condition $y\leq x.$ Analogously, if $A=[a_{i,j}]$ is poset matrix of size $n\times n,$ then it's dual poset matrix $A^{*}=[a_{{n+1-j},{n+1-i}}].$ This is reflected in partial composition operation on poset matrices as  summarized in the theorem that follows.

\begin{theorem}\label{dual}
Let $\text{A}_n$ and $\text{B}_m$ be poset matrices of order $n$ and $m$ respectively and let $i$  be any label of $A_n$. If $\text{A}_n\;\scalebox{0.7}{$\square$}_i\;\text{B}_m=\text{C}_{n+m-1}$ then for every $i,$  $\text{A}^{*}_n\;\scalebox{0.7}{$\square$}_{n-i+1}\;\text{B}_m=\text{C}^{*}_{n+m-1}$ where $\text{C}^{*}$ denotes the dual poset matrix of $\text{C}.$ 
\end{theorem}

\noindent Similar results relating poset matrix construction to their dual poset matrices can be formulated using the second and third partial composition operations of poset matrices.

\begin{theorem}\label{dualsecond}
Let $\text{A}_n$ and $\text{B}_m$ be poset matrices of order $n$ and $m$ respectively and let $i$  be any label of $A_n$. If $\text{A}_n\;\scalebox{0.7}{$\blacktriangle$}_i\;\text{B}_m=\text{C}_{n+m-1}$ then for every $i,$  $\text{A}^{*}_n\;\scalebox{0.7}{$\blacktriangledown$}_{n-i+1}\;\text{B}_m=\text{C}^{*}_{n+m-1}$ where $\text{C}^{*}$ denotes the dual poset matrix of $\text{C}.$ 
\end{theorem}

The two connected dual poset matrices of order $3$ are as follows.

\begin{example}
\end{example}

\begin{enumerate}
\item $\begin{blockarray}{ccc}
1 & 2 \\
\begin{block}{(cc)c}
 1 & 0  &1\\
 1 & 1 &2 \\
\end{block}
\end{blockarray}\scalebox{0.7}{$\square$}_1\;\;\begin{blockarray}{ccc}
1 & 2 \\
\begin{block}{(cc)c}
 1 & 0&1 \\
 0 & 1&2\\
\end{block}
\end{blockarray}=
\begin{blockarray}{cccc}
1&2& 3 \\
\begin{block}{(ccc)c}
 1 & 0 & 0&1 \\
 0 & 1 & 0&2 \\
 1 & 1 & 1&3 \\
\end{block}
\end{blockarray}\quad\Leftrightarrow\quad$
\begin{minipage}[c]{.20\textwidth}
\begin{tikzpicture}
  [scale=.7,auto=center,every node/.style={circle,fill=blue!20}]
\node(d1) at (10,1.7) {3};
\node(d6) at (11,0.5)  {2};
\node(d7) at (9,0.5) {1};
\draw (d6) -- (d1);
\draw (d1) -- (d7);

 \end{tikzpicture}

\end{minipage}

\item $\begin{blockarray}{ccc}
1 & 2 \\
\begin{block}{(cc)c}
 1 & 0  &1\\
 1 & 1 &2 \\
\end{block}
\end{blockarray}\scalebox{0.7}{$\square$}_2\;\;\begin{blockarray}{ccc}
1 & 2 \\
\begin{block}{(cc)c}
 1 & 0&1 \\
 0 & 1&2\\
\end{block}
\end{blockarray}=
\begin{blockarray}{cccc}
1&2 & 3 \\
\begin{block}{(ccc)c}
 1 & 0 & 0&1 \\
 1 & 1 & 0&2 \\
 1 & 0 & 1&3 \\
\end{block}
\end{blockarray}\quad\Leftrightarrow\quad$
\begin{minipage}[c]{.20\textwidth}
\begin{tikzpicture}
  [scale=.7,auto=center,every node/.style={circle,fill=blue!20}]
\node(d1) at (10,0.5) {1};
\node(d6) at (11,1.7)  {3};
\node(d7) at (9,1.7) {2};
\draw (d6) -- (d1);
\draw (d1) -- (d7);

 \end{tikzpicture}

\end{minipage}

\end{enumerate}

\begin{example}
\end{example}

The $10$ connected poset matrices of order $4$  and their duals are as follows.

\begin{enumerate}
\item $\begin{blockarray}{cccc}
1 & 2&3\\
\begin{block}{(ccc)c}
 1 & 0&0  &1\\
 1 & 1&0 &2\\
 1 & 1&1 &3 \\
\end{block}
\end{blockarray}\scalebox{0.7}{$\square$}_3\;\;
\begin{blockarray}{ccc}
1 & 2 \\
\begin{block}{(cc)c}
 1 & 0&1 \\
 0 & 1&2\\
\end{block}
\end{blockarray}=
\begin{blockarray}{ccccc}
1&2 & 3 &4\\
\begin{block}{(cccc)c}
 1 & 0 & 0&0 &1\\
 1 & 1 & 0&0 &2\\
 1 & 1 & 1&0 &3\\
 1 & 1 & 0&1 &4\\
\end{block}
\end{blockarray}\quad\Leftrightarrow\quad$
\begin{minipage}[c]{.20\textwidth}
\begin{tikzpicture}
  [scale=.7,auto=center,every node/.style={circle,fill=blue!20}]
\node(d1) at (10,4) {2};
\node(d6) at (11,5)  {4};
\node(d7) at (9,5) {3};
\node(d8) at (10,2.5) {1};
\draw (d6) -- (d1);
\draw (d1) -- (d7);
\draw (d1) -- (d8);

 \end{tikzpicture}

\end{minipage}

\item $\begin{blockarray}{cccc}
1 & 2&3\\
\begin{block}{(ccc)c}
 1 & 0&0  &1\\
 0 & 1&0 &2\\
 1 & 1&1 &3 \\
\end{block}
\end{blockarray}\scalebox{0.7}{$\square$}_3\;\;\begin{blockarray}{ccc}
1 &2 \\
\begin{block}{(cc)c}
 1 & 0&1 \\
 1 & 1&2\\
\end{block}
\end{blockarray}=
\begin{blockarray}{ccccc}
1&2 & 3 &4\\
\begin{block}{(cccc)c}
 1 & 0 & 0&0 &1\\
 0 & 1 & 0&0 &2\\
 1 & 1 & 1&0 &3\\
 1 & 1 & 1&1 &4\\
\end{block}
\end{blockarray}\quad\Leftrightarrow\quad$
\begin{minipage}[c]{.20\textwidth}
\begin{tikzpicture}
  [scale=.7,auto=center,every node/.style={circle,fill=blue!20}]

\node(d4) at (10,2) {4};
\node(d1) at (10,0.5) {3};
\node(d6) at (11,-0.3)  {2};
\node(d7) at (9,-0.3) {1};
\draw (d6) -- (d1);
\draw (d1) -- (d7);
\draw (d1) -- (d4);

 \end{tikzpicture}

\end{minipage}

\item $\begin{blockarray}{cccc}
1 & 2&3\\
\begin{block}{(ccc)c}
 1 & 0&0  &1\\
 1 & 1&0 &2\\
 1 & 1&1 &3 \\
\end{block}
\end{blockarray}\scalebox{0.7}{$\square$}_1\;\;
\begin{blockarray}{ccc}
1 & 2 \\
\begin{block}{(cc)c}
 1 & 0&1 \\
 1 & 1&2\\
\end{block}
\end{blockarray}=
\begin{blockarray}{ccccc}
1&2 & 3 &4\\
\begin{block}{(cccc)c}
 1 & 0 & 0&0 &1\\
 1 & 1 & 0&0 &2\\
 1 & 1 & 1&0 &3\\
 1 & 1 & 1&1 &4\\
\end{block}
\end{blockarray}\quad\Leftrightarrow\quad$
\begin{minipage}[c]{.20\textwidth}
\begin{tikzpicture}
  [scale=.7,auto=center,every node/.style={circle,fill=blue!20}]
\node(d1) at (10,4) {3};
\node(d6) at (10,3)  {2};
\node(d7) at (10,5) {4};
\node(d8) at (10,2) {1};
\draw (d6) -- (d1);
\draw (d1) -- (d7);
\draw (d6) -- (d8);

 \end{tikzpicture}

\end{minipage}

\item $
\begin{blockarray}{ccc}
1 & 2 \\
\begin{block}{(cc)c}
 1 & 0&1 \\
 1 & 1&2\\
\end{block}
\end{blockarray}
\scalebox{0.7}{$\square$}_1\;\;\begin{blockarray}{cccc}
1 & 2&3\\
\begin{block}{(ccc)c}
 1 & 0&0  &1\\
 1 & 1&0 &2\\
 1 & 0&1 &3 \\
\end{block}
\end{blockarray}
=\begin{blockarray}{ccccc}
1&2 & 3 &4\\
\begin{block}{(cccc)c}
 1 & 0 & 0&0 &1\\
 1 & 1 & 0&0 &2\\
 1 & 0 & 1&0 &3\\
 1 & 1 & 1&1 &4\\
\end{block}
\end{blockarray}\quad\Leftrightarrow\quad$
\begin{minipage}[c]{.20\textwidth}
\begin{tikzpicture}
  [scale=.7,auto=center,every node/.style={circle,fill=blue!20}]

\node[rectangle] (c1) at (6.2,3) {2};
  \node[rectangle] (c2) at (7,2) {1};
  \node[rectangle] (c4) at (7.8,3) {3};
  \node[rectangle] (c8) at (7,4)  {4};

      \draw (c1) -- (c2);
        \draw (c2) -- (c4);
    \draw (c8) -- (c4);
   \draw (c8) -- (c1);

 \end{tikzpicture}

\end{minipage}

\item $\begin{blockarray}{cccc}
1 & 2&3\\
\begin{block}{(ccc)c}
 1 & 0&0  &1\\
 1 & 1&0 &2\\
 1 & 0&1 &3 \\
\end{block}
\end{blockarray}\scalebox{0.7}{$\square$}_1\;\;
\begin{blockarray}{ccc}
1& 2 \\
\begin{block}{(cc)c}
 1 & 0&1 \\
 0 & 1&2\\
\end{block}
\end{blockarray}=
\begin{blockarray}{ccccc}
1&2 & 3 &4\\
\begin{block}{(cccc)c}
 1 & 0 & 0&0 &1\\
 0 & 1 & 0&0 &2\\
 1 & 1 & 1&0 &3\\
 1 & 1 & 0&1 &4\\
\end{block}
\end{blockarray}\quad\Leftrightarrow\quad$
\begin{minipage}[c]{.20\textwidth}
\begin{tikzpicture}
  [scale=.7,auto=center,every node/.style={circle,fill=blue!20}]
\node(d1) at (10,5) {3};
\node(d8) at (10,3) {1};
\node(d4) at (12,5) {4};
\node(d5) at (12,3) {2};
\draw (d4) -- (d5);
\draw (d1) -- (d8);
\draw (d1) -- (d5);
\draw (d4) -- (d8);

 \end{tikzpicture}

\end{minipage}

\item $\begin{blockarray}{cccc}
1 & 2&3\\
\begin{block}{(ccc)c}
 1 & 0&0  &1\\
 1 & 1&0 &2\\
 1 & 0&1 &3 \\
\end{block}
\end{blockarray}\scalebox{0.7}{$\square$}_2\;\;\begin{blockarray}{ccc}
1 &2 \\
\begin{block}{(cc)c}
 1 & 0&1 \\
 0 & 1&2\\
\end{block}
\end{blockarray}
=\begin{blockarray}{ccccc}
1&2 & 3 &4\\
\begin{block}{(cccc)c}
 1 & 0 & 0&0 &1\\
 1 & 1 & 0&0 &2\\
 1 & 0 & 1&0 &3\\
 1 & 0 & 0&1 &4\\
\end{block}
\end{blockarray}\quad\Leftrightarrow\quad$
\begin{minipage}[c]{.20\textwidth}
\begin{tikzpicture}
  [scale=.7,auto=center,every node/.style={circle,fill=blue!20}]

\node[rectangle] (c1) at (5.5,3) {2};
  \node[rectangle] (c2) at (7,1) {1};
  \node[rectangle] (c4) at (7,3) {3};
   \node[rectangle] (c5) at (8.5,3) {4};
      \draw (c5) -- (c2);
      \draw (c1) -- (c2);
        \draw (c2) -- (c4);

 \end{tikzpicture}

\end{minipage}

\item $\begin{blockarray}{cccc}
1 & 2&3\\
\begin{block}{(ccc)c}
 1 & 0&0  &1\\
 0 & 1&0 &2\\
 1 & 1&1 &3 \\
\end{block}
\end{blockarray}\scalebox{0.7}{$\square$}_1\;\;
\begin{blockarray}{ccc}
1 & 2 \\
\begin{block}{(cc)c}
 1 & 0&1 \\
 0 & 1&2\\
\end{block}
\end{blockarray}=
\begin{blockarray}{ccccc}
1&2 & 3 &4\\
\begin{block}{(cccc)c}
 1 & 0 & 0&0 &1\\
 0 & 1 & 0&0 &2\\
 0 & 0 & 1&0 &3\\
 1 & 1 & 1&1 &4\\
\end{block}
\end{blockarray}\quad\Leftrightarrow\quad$
\begin{minipage}[c]{.20\textwidth}
\begin{tikzpicture}
  [scale=.7,auto=center,every node/.style={circle,fill=blue!20}]

\node(d4) at (10,-0.3) {2};
\node(d1) at (10,1) {4};
\node(d6) at (11,-0.3)  {1};
\node(d7) at (9,-0.3) {3};
\draw (d6) -- (d1);
\draw (d1) -- (d7);
\draw (d1) -- (d4);

 \end{tikzpicture}

\end{minipage}

\item $\begin{blockarray}{cccc}
1 & 2&3\\
\begin{block}{(ccc)c}
 1 & 0&0  &1\\
 0 & 1&0 &2\\
 1 & 1&1 &3 \\
\end{block}
\end{blockarray}\scalebox{0.7}{$\square$}_1\;\;\begin{blockarray}{ccc}
1 &2 \\
\begin{block}{(cc)c}
 1 & 0&1 \\
 1 & 1&2\\
\end{block}
\end{blockarray}
=\begin{blockarray}{ccccc}
1&2 & 3 &4\\
\begin{block}{(cccc)c}
 1 & 0 & 0&0 &1\\
 1 & 1 & 0&0 &2\\
 0 & 0 & 1&0 &3\\
 1 & 1 & 1&1 &4\\
\end{block}
\end{blockarray}\quad\Leftrightarrow\quad$
\begin{minipage}[c]{.20\textwidth}
\begin{tikzpicture}
  [scale=.7,auto=center,every node/.style={circle,fill=blue!20}]

\node[rectangle] (c1) at (6.2,3) {3};
  \node[rectangle] (c2) at (7.8,1.7) {1};
  \node[rectangle] (c4) at (7.8,3) {2};
  \node[rectangle] (c8) at (7,4)  {4};

      \draw (c4) -- (c2);
    \draw (c8) -- (c4);
   \draw (c8) -- (c1);

 \end{tikzpicture}

\end{minipage}

\item $\begin{blockarray}{cccc}
1 & 2&3\\
\begin{block}{(ccc)c}
 1 & 0&0  &1\\
 1 & 1&0 &2\\
 1 & 0&1 &3 \\
\end{block}
\end{blockarray}\scalebox{0.7}{$\square$}_2\;\;\begin{blockarray}{ccc}
1 & 2 \\
\begin{block}{(cc)c}
 1 & 0&1 \\
 1 & 1&2\\
\end{block}
\end{blockarray}
=\begin{blockarray}{ccccc}
1&2 & 3 &4\\
\begin{block}{(cccc)c}
 1 & 0 & 0&0 &1\\
 1 & 1 & 0&0 &2\\
 1 & 1 & 1&0 &3\\
 1 & 0 & 0&1 &4\\
\end{block}
\end{blockarray}\quad\Leftrightarrow\quad$
\begin{minipage}[c]{.20\textwidth}
\begin{tikzpicture}
  [scale=.7,auto=center,every node/.style={circle,fill=blue!20}]

\node[rectangle] (c1) at (6.2,2.5) {2};
  \node[rectangle] (c2) at (7,1.5) {1};
  \node[rectangle] (c4) at (7.8,2.5) {4};
  \node[rectangle] (c8) at (6.2,3.7)  {3};

      \draw (c1) -- (c2);
       \draw (c4) -- (c2);
   \draw (c8) -- (c1);

 \end{tikzpicture}

\end{minipage}

\item $\begin{blockarray}{cccc}
1 & 2&3\\
\begin{block}{(ccc)c}
 1 & 0&0  &1\\
 1 & 1&0 &2\\
 1 & 0&1 &3 \\
\end{block}
\end{blockarray}\blacktriangle_2\;\;
\begin{blockarray}{ccc}
1 & 2 \\
\begin{block}{(cc)c}
 1 & 0&1 \\
 1 & 1&2\\
\end{block}
\end{blockarray}=
\begin{blockarray}{ccccc}
1&2 & 3 &4\\
\begin{block}{(cccc)c}
 1 & 0 & 0&0 &1\\
 0 & 1 & 0&0 &2\\
 1 & 1 & 1&0 &3\\
 1 & 0 & 0&1 &4\\
\end{block}
\end{blockarray}\quad\Leftrightarrow\quad$
\begin{minipage}[c]{.20\textwidth}
\begin{tikzpicture}
  [scale=.7,auto=center,every node/.style={circle,fill=blue!20}]

\node(d1) at (9,1) {2};
\node(d6) at (9,3)  {3};
\node(d3) at (11,1)  {1};
\node(d5) at (11,3)  {4};
\draw (d6) -- (d1);
\draw (d3) -- (d5);
\draw (d1) -- (d5);

 \end{tikzpicture}

\end{minipage}

\end{enumerate}

Similarly, the $6$ disconnected poset matrices of order $4$ are generated as follows.

\begin{enumerate}

\item $\begin{blockarray}{cccc}
1 & 2&3\\
\begin{block}{(ccc)c}
 1 & 0&0& 1\\
 1 & 1&0 &2\\
  0 & 0&1&3\\
\end{block}
\end{blockarray}\scalebox{0.7}{$\square$}_3\;\;
\begin{blockarray}{ccc}
1 & 2 \\
\begin{block}{(cc)c}
1 & 0&1\\
1 & 1&2\\
\end{block}
\end{blockarray}=
\begin{blockarray}{ccccc}
1&2 & 3 &4\\
\begin{block}{(cccc)c}
 1 & 0 & 0&0&1 \\
 1 & 1 & 0&0&2 \\
 0 & 0 & 1&0&3 \\
 0 & 0 & 1&1&4 \\
\end{block}
\end{blockarray}\quad\Leftrightarrow\quad$
\begin{minipage}[c]{.20\textwidth}
 \begin{tikzpicture}
\node [rectangle] (bb1) at (3,5) {2};
\node [rectangle] (bb2) at (3,3)  {1};
\node [rectangle]  (bb3) at (4,5) {4};
\node [rectangle] (bb4) at (4,3)  {3};
\node [use](aa17) at (5,4) {$.$};
\draw (bb2) -- (bb1);
\draw (bb3) -- (bb4);
\end{tikzpicture}

\end{minipage}

\item $\begin{blockarray}{cccc}
1 & 2&3\\
\begin{block}{(ccc)c}
1 & 0&0&1  \\
 0 & 1&0&2 \\
 0 & 0&1&3\\
\end{block}
\end{blockarray}\scalebox{0.7}{$\square$}_1\;\;
\begin{blockarray}{ccc}
1 & 2 \\
\begin{block}{(cc)c}
1 & 0&1\\
0 & 1&2\\
\end{block}
\end{blockarray}=
\begin{blockarray}{ccccc}
1&2 & 3 &4\\
\begin{block}{(cccc)c}
  1 & 0 & 0&0&1 \\
 0 & 1 & 0&0&2 \\
 0 & 0 & 1&0&3 \\
 0 & 0 & 0&1&4 \\
\end{block}
\end{blockarray}\quad\Leftrightarrow\quad$
\begin{minipage}[c]{.20\textwidth}
\begin{tikzpicture}
  [scale=.7,auto=center,every node/.style={circle,fill=blue!20}]
\node(d1) at (9,1) {1};
\node(d6) at (10,1)  {2};
\node(d7) at (11,1) {3};
\node(d8) at (12,1) {4};
\end{tikzpicture}
 \end{minipage}

\item $\begin{blockarray}{cccc}
1 & 2&3\\
\begin{block}{(ccc)c}
 1 & 0&0&1  \\
 0 & 1&0&2\\
   0 & 1&1&3\\
\end{block}
\end{blockarray}\scalebox{0.7}{$\square$}_1\;\;
\begin{blockarray}{ccc}
1 & 2 \\
\begin{block}{(cc)c}
1 & 0& 1\\
0 & 1& 2\\
\end{block}
\end{blockarray}=
\begin{blockarray}{ccccc}
1&2 & 3 &4\\
\begin{block}{(cccc)c}
 1 & 0 & 0&0&1 \\
 0 & 1 & 0&0&2 \\
 0 & 0 & 1&0&3 \\
 0 & 0 & 1&1&4 \\
\end{block}
\end{blockarray}\quad\Leftrightarrow\quad$
\begin{minipage}[c]{.20\textwidth}
\begin{tikzpicture}
  [scale=.7,auto=center,every node/.style={circle,fill=blue!20}]
\node(d1) at (9,1) {3};
\node(d6) at (9,2.5)  {4};
\node(d7) at (10,1) {1};
\node(d7) at (11,1) {2};
\draw (d6) -- (d1);
\end{tikzpicture}
 \end{minipage}

\item $\begin{blockarray}{cccc}
1 & 2&3\\
\begin{block}{(ccc)c}
  1 & 0&0&1 \\
 1 & 1&0&2 \\
  0 & 0&1&3\\
\end{block}
\end{blockarray}\scalebox{0.7}{$\square$}_1\;\;
\begin{blockarray}{ccc}
1 & 2 \\
\begin{block}{(cc)c}
1 & 0&1\\
1 & 1&2\\
\end{block}
\end{blockarray}=
\begin{blockarray}{ccccc}
1&2 & 3 &4\\
\begin{block}{(cccc)c}
  1 & 0 & 0&0 &1\\
 1 & 1 & 0&0&2 \\
 1 & 1 & 1&0&3 \\
 0 & 0 & 0&1&4 \\
\end{block}
\end{blockarray}\quad\Leftrightarrow\quad$
\begin{minipage}[c]{.20\textwidth}
\begin{tikzpicture}
  [scale=.7,auto=center,every node/.style={circle,fill=blue!20}]
\node(d1) at (9,1) {1};
\node(d6) at (9,2.5)  {2};
\node(d8) at (9,4)  {3};
\node(d7) at (10,1) {4};
\draw (d6) -- (d1);
\draw (d6) -- (d8);
\end{tikzpicture}
 \end{minipage}

\item $\begin{blockarray}{cccc}
1 & 2&3\\
\begin{block}{(ccc)c}
  1 & 0&0&1  \\
 0 & 1&0&2 \\
  0 & 1&1&3\\
\end{block}
\end{blockarray}\scalebox{0.7}{$\square$}_2\;\;
\begin{blockarray}{ccc}
1 & 2 \\
\begin{block}{(cc)c}
1 & 0&1\\
0 & 1&2\\
\end{block}
\end{blockarray}=
\begin{blockarray}{ccccc}
1&2 & 3 &4\\
\begin{block}{(cccc)c}
   1 & 0 & 0&0&1 \\
 0 & 1 & 0&0&2 \\
 0 & 0 & 1&0&3 \\
 0 & 1 & 1&1&4 \\
\end{block}
\end{blockarray}\quad\Leftrightarrow\quad$
\begin{minipage}[c]{.20\textwidth}
\begin{tikzpicture}
  [scale=.7,auto=center,every node/.style={circle,fill=blue!20}]
\node(d1) at (10,1.7) {4};
\node(d6) at (11,0.5)  {3};
\node(d7) at (9,0.5) {2};
\node(d88) at (12,0.5) {1};
\draw (d6) -- (d1);
\draw (d1) -- (d7);

 \end{tikzpicture}
  \end{minipage}

\item $\begin{blockarray}{cccc}
1 & 2&3\\
\begin{block}{(ccc)c}
  1 & 0&0&1  \\
 0 & 1&0&2 \\
   0 & 1&1&3\\
\end{block}
\end{blockarray}\scalebox{0.7}{$\square$}_3\;\;
\begin{blockarray}{ccc}
1 & 2 \\
\begin{block}{(cc)c}
1 & 0&1\\
 0 & 1&2\\
\end{block}
\end{blockarray}=
\begin{blockarray}{ccccc}
1&2 & 3 &4\\
\begin{block}{(cccc)c}
     1 & 0 & 0&0 &1\\
 0 & 1 & 0&0 &2\\
 0 & 1 & 1&0 &3\\
 0 & 1 & 0&1&4 \\
\end{block}
\end{blockarray}\quad\Leftrightarrow\quad$
\begin{minipage}[c]{.20\textwidth}
\begin{tikzpicture}
  [scale=.7,auto=center,every node/.style={circle,fill=blue!20}]
\node(d1) at (10,0.5) {2};
\node(d6) at (11,1.7)  {4};
\node(d7) at (9,1.7) {3};
\node(d88) at (12,0.5) {1};
\draw (d6) -- (d1);
\draw (d1) -- (d7);

 \end{tikzpicture}
   \end{minipage}

\end{enumerate}

\begin{remark} It is worth noting that the poset matrices presented above can also be generated using the second and third partial composition operation on poset matrices. 
\end{remark}

\begin{definition}{\rm 
Let  $\text{A}=[a_{i,j}]$ be a poset matrix of size $n\times n$. If $a_{i,j}=1$ for all $i\in[n]$ whenever $i\leq j,$  then we call $\text{A}$ a totally connnected poset matrix.}
\end{definition}
\begin{definition}{\rm 
Let  $\text{A}=[a_{i,j}]$ be a poset matrix of size $n\times n.$ If $a_{i,j}=0$ for all $i\in[n]$ whenever $i\neq j,$ Then we call $\text{A}$ a totally disconnected poset matrix.}
\end{definition}

%the j^{th} row of size 1x(k-1) and the  j^{th}  column of size (n-k+1)x1

\begin{remark}{\rm 
For the purpose of clarity on the notations used in the theorems that follow, given a poset matrix $A$ of size $n\times n,$ we refer to $A_{(j)}[\{k,...,n\}\mid\{1,..,k-1\}]$ and $A^{(j)}[\{k,...,n\}\mid\{1,..,k-1\}]$ to denote the $j^{th}$ row of size $1\times(k-1)$ and the $j^{th}$ column of size $(n-k+1)\times 1$ respectively of the submatrix $A[\{k,...,n\}\mid\{1,..,k-1\}]$ derived from the matrix $A.$ If $k=2$ and $n=5$ then the set $\{k,...,n\}$ is equivalent to $\{2,3,4,5\}.$} We denote by $\mathcal{PM}(n)$ the set of all $n\times n$ poset matrices and let $$\mathcal{PM}=\bigsqcup\limits_{n\ge 1}\mathcal{PM}(n).$$

\end{remark}

\begin{theorem}\label{connectedsubposet} Let $A$ be a connected poset matrix of size $n\times n$ and let $\text{B}$ be a totally connected poset matrix of size $m\times m.$ Then the following holds.
\begin{enumerate}
\item  If there exist a totally connected subposet matrix $A[\alpha]$ defined on $\alpha=\{1,...,k\}$ with $1<k<n$ such that the submatrix $D:=A[\{k+1,...,n\}|\{1,...,k\}]$ satisfies the condition that $D^{(p)}=D^{(p+1)}$ for each $p\in\{1,...,k-1\}.$ Then $\text{A}\square_i \text{B}=\text{A}\square_r \text{B}$ whenever $i,r\in\alpha.$
\item If there exist a totally connected subposet matrix $A[\alpha]$ defined on $\alpha=\{k,...,n\}$ with $k\geq 2$ such that the submatrix $D:=A[\{k,...,n\}|\{1,...,k-1\}]$ satisfies the condition that $D_{(p)}=D_{(p+1)}$ for each $p\in\{1,...,n-k\}.$ Then $\text{A}\square_i \text{B}=\text{A}\square_r \text{B}$ whenever $i,r\in\alpha.$
\item   If there exist a totally connected subposet matrix $A[\alpha]$ defined on $\alpha=\{d,...,k\}$ with  $1<d<k<n$ such that the submatrix $D:=A[\{k,...,n\}|\{1,...,k-1\}]$ satisfies the condition $D_{(p)}=D_{(p+1)}$ for each $p\in\{1,...,n-k\}$ where $D_{(p)}$ is a   \textbf{1}-vector or a \textbf{0}-vector of size $1 \times(k-1).$ Then $\text{A}\square_i \text{B}=\text{A}\square_r \text{B}$ whenever $i,r\in\alpha.$
\end{enumerate}
\end{theorem}
\begin{proof}
Let $j,j+1\in\alpha.$ In (1), it suffices to show that the output poset matrix $A\square_{j} B=A\square_{j +1}B$ for $j\in\{1,...,k\}.$   Consider that the output poset matrix structure of $A\square_{j}B=E$ can be partitioned into three submatrices $E1, E2$ and $E3$ such that:
\begin{itemize}
\item $E1=E[\{1,...,k+m-1\}]$  of size $(k+m-1)\times (k+m-1).$
\item $E2=E[\{k+m,...,n+m-1\}\mid\{1,...,k+m-1\}]$ of size $(n-k)\times (k+m-1).$
\item $E3=E[\{k+m,...,n+m-1\}]$ of size $(n-k)\times(n-k).$
\end{itemize}
$E1$ can be subpartitioned such that $E[\{j,...,j+m\}]=B$, $E[\{j,...,j+m\}\mid\{1,...,j-1\} ]$, $E[\{j+m+1,....,k+m-1\}\mid\{1,...,k+m-1\}]$ and $E[\{1,...,j-1\}]$  have entries of all $1$'s that lie on and below the main diagonal of $E$ since these entries are derived from $A[\{1,...,k\}]$ which is a totally connected subposet matrix.

\noindent $E2$ can be subpartitioned such that $E[\{k+m,...,n+m-1\}\mid\{1,...,j-1\}]=A[\{k+1,...,n+m-1\}\mid\{1,...,j-1\}], E[\{k+m,...,n+m-1\}\mid\{j,...,j+m\}]={\mathbbm{1}}_{m}\otimes A^{(j)}[\{k+1,..n\}\mid\{1,..,k\}]$, $ E[\{k+m,...,n+m-1\}\mid\{j+m+1,...,k\}]=A[\{k+1,...,n+m-1\}\mid\{j+m+1,...,k]$ and $E3=A[\{k+1,...,n+m-1\}].$

In the adjacent insertion point $j+1$, similarly consider that the output poset matrix structure of $A\square_{j+1}B=F$ can be partitioned into three submatrices $F1, F2$ and $F3$ such that:
\begin{itemize}
\item $F1=F[\{1,...,k+m-1\}]$  of size $(k+m-1)\times (k+m-1).$
\item $F2=F[\{k+m,...,n+m-1\}\mid\{1,...,k+m-1\}]$ of size $(n-k)\times (k+m-1).$
\item $F3=F[\{k+m,...,n+m-1\}]$ of size $(n-k)\times(n-k).$
\end{itemize}
$F1$ can be subpartitioned such that $F[\{j+1,...,j+1+m\}]=B$, $F[\{j+1,...,j+1+m\}\mid\{1,...,j\} ]$, $F[\{j+m+2,....,k+m-1\}\mid\{1,...,k+m-1\}]$ and $F[\{1,...,j\}]$  have entries of all $1$'s that lie on and below the main diagonal of $F$ since these entries are derived from $A[\{1,...,k\}]$ which is a totally connected subposet matrix.

\noindent $F2$ can be subpartitioned such that $F[\{k+m,...,n+m-1\}\mid\{1,...,j\}]=A[\{k+1,...,n+m-1\}\mid\{1,...,j\}], F[\{k+m,...,n+m-1\}\mid\{j+1,...,j+1+m\}]={\mathbbm{1}}_{m}\otimes A^{(j+1)}[\{k+1,..n\}\mid\{1,..,k\}]$, $ F[\{k+m,...,n+m-1\}\mid\{j+m+2,...,k\}]=A[\{k+1,...,n+m-1\}\mid\{j+m+2,...,k]$ and $F3=A[\{k+1,...,n+m-1\}].$

From the matrix structure of $E$ and $F$ presented above, we get the following.

\begin{itemize}
\item $E1=F1$ since all the subpartitions of $E1$ and $F1$ are equal.
\item $E2=F2$ since as the subpartitions can be equal whenever ${\mathbbm{1}}_{m}\otimes A^{(j)}[\{k+1,..n\}\mid\{1,..,k\}]={\mathbbm{1}}_{m}\otimes A^{(j+1)}[\{k+1,..n\}\mid\{1,..,k\}].$ This condition holds when each of the columns in the submatrix $A[\{k+1,..n\}\mid\{1,..,k\}]$ have their correponding entries equal.
\item $E3=F3.$
\end{itemize}
 Thus,  $A\square_{j} B=A\square_{j +1}B$ for $j\in\{1,...,k\}.$

Let $j,j+1\in\alpha.$ In (2) it suffices to show that the output poset matrix $A\square_{j} B=A\square_{j +1}B$ for $j\in\{k,...,n\}.$  Consider that the output poset matrix structure of $A\square_{j}B=G$ can be partitioned into three submatrices $G1, G2$ and $G3$ such that:

\begin{itemize}
\item $G1=G[\{1,...,k-1\}]$ of size $(k-1)\times(k-1).$
\item  $G2=G[\{k,...,n+m-1\}\mid\{1,...,k-1\}]$ of size $(n+m-k)\times(k-1).$
\item $G3=G[\{k,...,n+m-1\}]$ of size $(n+m-k)\times(n+m-k)$.
\end{itemize}
We note that $G1=A[\{1,...,k-1\}].$ $G2$ can be subpartitioned such that  $G[\{j,...,j+m\}\mid\{1,...,k-1\}]={\mathbbm{1}}_m^T\otimes A_{(j)}[\{k,...,n\}\mid\{1,...,k-1\}]$ and $G[\{j+m+1,...,n+m-1\}\mid\{1,...,k-1\}]=A[\{j+1,...,n\}\mid\{1,...,k-1].$ $G3$ can be subpartitioned such that $G[\{j,...j+m\}]=B,$ $G[\{k,..,j-1\}]$, and $G[\{j+m+1,...,n+m-1\}]$ have entries of all $1$'s that lie on and below the main diagonal of $G$ since these entries are derived from $A[\{k,...,n\}]$ which is a totally connected subposet matrix.

In the adjacent insertion point $j+1$, similarly consider that the output poset matrix structure of $A\square_{j+1}B=H$ can be partitioned into three submatrices $H1, H2$ and $H3$ such that:
\begin{itemize}
\item $H1=H[\{1,...,k-1\}]$ of size $(k-1)\times(k-1).$
\item  $H2=H[\{k,...,n+m-1\}\mid\{1,...,k-1\}]$ of size $(n+m-k)\times(k-1).$
\item $H3=H[\{k,...,n+m-1\}]$ of size $(n+m-k)\times(n+m-k)$.
\end{itemize}
Based on the matrix structures for $G$ and $H$ we obtain the following.
\begin{itemize}
\item $H1=A[\{1,...,k-1\}]=G1.$
\item By replacing $j:=j+1$ in $G2$ and noting that ${\mathbbm{1}}_m^T\otimes A_{(j)}[\{k,...,n\}\mid\{1,...,k-1\}]={\mathbbm{1}}_m^T\otimes A_{(j+1)}[\{k,...,n\}\mid\{1,...,k-1\}]$ is satisfied  when each of the rows in the submatrix $A[\{k,..,n\}\mid\{1,..,k-1\}]$ have their correponding entries equal. Therefore whenever this condition holds $G2=H2.$
\item $G3=H3$ since as replacing $j$ in $E3$ with $j+1,$ preserves the same matrix structure in $H3.$
\end{itemize}
 Thus,  $A\square_{j} B=A\square_{j +1}B$ for each $j\in\{k,...,n\}.$

\noindent The proof of (3) follows from the arguments in (1) and (2).

\end{proof}

\begin{example} {\rm Consider the poset matrices}
\end{example}
 $$\text{A}=\begin{blockarray}{ccccc}
1 & 2&3&4\\
\begin{block}{(cccc)c}
 1 & 0&0  &0&1\\
 1 & 1&0 &0&2\\
 1 & 1&1 &0&3 \\
 1&1&0&1&4\\
\end{block}
\end{blockarray}\qquad\text{B}=\begin{blockarray}{ccc}
1 & 2 \\
\begin{block}{(cc)c}
 1 & 0&1 \\
 1 & 1&2\\
\end{block}
\end{blockarray}\qquad\text{C}=\begin{blockarray}{cccccc}
1&2 & 3 &4&5\\
\begin{block}{(ccccc)c}
 1 & 0 & 0&0&0 &1\\
 1 & 1 & 0&0& 0&2\\
 1 & 1 & 1&0 &0&3\\
 1&1&1&1&0&4\\
 1 & 1 & 1&0&1 &5\\
\end{block}
\end{blockarray}$$

It can be verified that $\text{A}\square_1 \text{B}=\text{A}\square_2\text{B}=\text{C}.$ On the other hand, $\text{A}\square_1 \text{B}\neq\text{A}\square_3 \text{B}\neq\text{A}\square_4 \text{B}$
\begin{remark}{\rm
The subposet matrix $A[\alpha]$ with  $\alpha=\{1,2\}$ forms a totally connected subposet matrix of $A.$ The submatrix  $D=A[{3,4}|{1,2}]$ has $2$ equal columns. By Theorem \ref{connectedsubposet} the output poset matrices from insertion at the labels $\{1,2\}$ are identical. On the other hand, consider the totally connected subposet matrix  $A[\alpha]$ with  $\alpha=\{1,2,3\}.$ In this case its associated submatrix $D=A[\{4\}|\{1,2,3\}]$ do not have all equal columns since the entry of the third column of $D$ is different from the first and second column. By Theorem \ref{connectedsubposet} its output poset matrix from the square partial composition operation at insertion point $3$.would be different from the output poset matrices at insertion points $1$ and $2$ of the input poset matrix $A.$}

\end{remark}

\begin{example}{\rm Consider the poset matrices}
\end{example}
 $$\text{A}=\begin{blockarray}{ccccc}
1 & 2&3&4\\
\begin{block}{(cccc)c}
 1 & 0&0  &0&1\\
 1 & 1&0 &0&2\\
 0 & 0&1 &0&3 \\
 1&1&1&1&4\\
\end{block}
\end{blockarray}\qquad\text{B}=\begin{blockarray}{ccc}
1 & 2 \\
\begin{block}{(cc)c}
 1 & 0&1 \\
 1 & 1&2\\
\end{block}
\end{blockarray}\qquad\text{C}=\begin{blockarray}{cccccc}
1&2 & 3 &4&5\\
\begin{block}{(ccccc)c}
 1 & 0 & 0&0&0 &1\\
 1 & 1 & 0&0& 0&2\\
 1 & 1 & 1&0 &0&3\\
 0&0&0&1&0&4\\
 1 & 1 & 1&1&1 &5\\
\end{block}
\end{blockarray}$$

It can be verified that $\text{A}\square_1 \text{B}=\text{A}\square_2\text{B}=\text{C}.$ On the other hand, $\text{A}\square_1 \text{B}\neq\text{A}\square_3 \text{B}\neq\text{A}\square_4 \text{B}$
\begin{remark}
{\rm The totally connected subposet matrix $A[\alpha]$ with  $\alpha=\{1,2\}$ forms a totally connected subposet matrix of $A.$ The submatrix  $D=A[{3,4}|{1,2}]$ has $2$ equal coluns. By Theorem \ref{connectedsubposet}(i) the outposet matrix from insertion at the labels $\{1,2\}$ are identical. On the other hand,  the totally connected subposet matrix $A[\alpha]$ with  $\alpha=\{3,4\}$ fdoes not result in identical output poset matrices at insertion points $3$ and $4.$  In this case, it can be observed that the associated submatrix $D=A[{3,4}|{1,2}]$ does not have equal rows as required by Theorem \ref{connectedsubposet}(ii).}

\end{remark}

\begin{example}{\rm Consider the poset matrices}
\end{example}
 $$\text{A}=\begin{blockarray}{ccccc}
1 & 2&3&4\\
\begin{block}{(cccc)c}
 1 & 0&0  &0&1\\
 1 & 1&0 &0&2\\
 1 & 0&1 &0&3 \\
 1&0&1&1&4\\
\end{block}
\end{blockarray}\qquad\text{B}=\begin{blockarray}{ccc}
1 & 2 \\
\begin{block}{(cc)c}
 1 & 0&1 \\
 1 & 1&2\\
\end{block}
\end{blockarray}\qquad\text{C}=\begin{blockarray}{cccccc}
1&2 & 3 &4&5\\
\begin{block}{(ccccc)c}
 1 & 0 & 0&0&0 &1\\
 1 & 1 & 0&0& 0&2\\
 1 & 0 & 1&0 &0&3\\
 1&0&1&1&0&4\\
 1 & 0 & 1&1&1 &5\\
\end{block}
\end{blockarray}$$

It can be verified that $\text{A}\square_3 \text{B}=\text{A}\square_4\text{B}=\text{C}.$ On the other hand, $\text{A}\square_3 \text{B}\neq\text{A}\square_1 \text{B}\neq\text{A}\square_2 \text{B}.$
\begin{remark}{\rm The totally connected subposet matrix $A[\alpha]$ with  $\alpha=\{3,4\}$ forms a totally connected subposet matrix of $A.$ The submatrix  $D=A[{3,4}|{1,2}]$ has $2$ equal rows. By Theorem \ref{connectedsubposet}(ii) the outposet matrix from insertion at the labels $\{3,4\}$ are identical. On the other hand,  the totally connected subposet matrix $A[\alpha]$ with  $\alpha=\{1,2\}$ does not result in identical output poset matrices at insertion points $1$ and $2.$  In this case, it can be observed that the associated submatrix $D=A[{3,4}|{1,2}]$ does not have equal columns as required by Theorem \ref{connectedsubposet}(i).
}

\end{remark}

\begin{theorem}\label{disconnectedsubposet} Let $A\in {\cal PM}(n)$ and let $\text{B}$ be a totally disconnected poset matrix of size $m\times m.$  Then the following holds.
\begin{enumerate}
\item  If there exist a totally disconnected subposet matrix $A[\alpha]$ defined on $\alpha=\{1,...,k\}$ with $1<k<n$ such that the submatrix $D:=A[\{k+1,...,n\}|\{1,...,k\}]$ satisfies the condition $D^{(p)}=D^{(p+1)}$ for each $p\in\{1,...,k-1\}.$ Then $\text{A}\square_i \text{B}=\text{A}\square_r \text{B}$ whenever $i,r\in\alpha.$
\item If there exist a totally disconnected subposet matrix $A[\alpha]$ defined on $\alpha=\{k,...,n\}$ with $k\geq 2$ such that the submatrix $D:=A[\{k,...,n\}|\{1,...,k-1\}]$ satisfies the condition $D_{(p)}=D_{(p+1)}$ for each $p\in\{1,...,n-k\}.$ Then $\text{A}\square_i \text{B}=\text{A}\square_r \text{B}$ whenever $i,r\in\alpha.$
\item   If there exist a totally disconnected subposet matrix $A[\alpha]$ defined on $\alpha=\{d,...,k\}$ with $1<d<k<n$ such that the submatrix $D:=A[\{k,...,n\}|\{1,...,k-1\}]$ satisfies the condition $D_{(p)}=D_{(p+1)}$ for each $p\in\{1,...,n-k\}$ where $D_{(p)}$ is a   \textbf{1}-vector or a \textbf{0}-vector of size $1 \times(k-1).$  Then $\text{A}\square_i \text{B}=\text{A}\square_r \text{B}$ whenever $i,r\in\alpha.$
\end{enumerate}
\end{theorem}

\begin{proof}
Similar to Theorem \ref{connectedsubposet}.
\end{proof}

\begin{example}
\end{example}

Consider the poset matrices $$\text{A}=\begin{blockarray}{ccccc}
1 & 2&3&4\\
\begin{block}{(cccc)c}
 1 & 0&0  &0&1\\
 1 & 1&0 &0&2\\
 1 & 1&1 &0&3 \\
 1&1&0&1&4\\
\end{block}
\end{blockarray}\qquad\text{B}=\begin{blockarray}{ccc}
1 & 2 \\
\begin{block}{(cc)c}
 1 & 0&1 \\
 0 & 1&2\\
\end{block}
\end{blockarray}\qquad\text{C}=\begin{blockarray}{cccccc}
1&2 & 3 &4&5\\
\begin{block}{(ccccc)c}
 1 & 0 & 0&0&0 &1\\
 1 & 1 & 0&0& 0&2\\
 1 & 1 & 1&0 &0&3\\
 1&1&0&1&0&4\\
 1 & 1 & 0&0&1 &5\\
\end{block}
\end{blockarray}$$

It can be verified that $\text{A}\square_3 \text{B}=\text{A}\square_4\text{B}=\text{C}.$ On the other hand, $\text{A}\square_1 \text{B}\neq\text{A}\square_2\text{B}\neq\text{A}\square_4 \text{B}.$

\begin{remark}
{\rm The subposet matrix $A[\alpha]$ with  $\alpha=\{3,4\}$ forms a totally disconnected subposet matrix of $A.$ The submatrix  $D=A[{3,4}|{1,2}]$ has $2$ equal rows. By Theorem \ref{disconnectedsubposet}(ii) the output poset matrix derived from the square PCO insertion at the labels $\{3,4\}$ are identical.}
\end{remark}

\begin{example}{\rm 
Consider the poset matrices}
\end{example}
 $$\text{A}=\begin{blockarray}{ccccc}
1 & 2&3&4\\
\begin{block}{(cccc)c}
 1 & 0&0  &0&1\\
 0 & 1&0 &0&2\\
 0 & 0&1 &0&3 \\
 1&1&1&1&4\\
\end{block}
\end{blockarray}\qquad\text{B}=\begin{blockarray}{ccc}
1 & 2 \\
\begin{block}{(cc)c}
 1 & 0&1 \\
 0 & 1&2\\
\end{block}
\end{blockarray}\qquad\text{C}=\begin{blockarray}{cccccc}
1&2 & 3 &4&5\\
\begin{block}{(ccccc)c}
 1 & 0 & 0&0&0 &1\\
 0 & 1 & 0&0& 0&2\\
 0 & 0 & 1&0 &0&3\\
 0&0&0&1&0&4\\
 1 & 1 & 1&1&1 &5\\
\end{block}
\end{blockarray}$$

It can be verified that $\text{A}\square_1 \text{B}=\text{A}\square_2\text{B}=\text{A}\square_3\text{B}=\text{C}.$ On the other hand, $\text{A}\square_1 \text{B}\neq\text{A}\square_4 \text{B}.$

\begin{remark}
{\rm The subposet matrix $A[\alpha]$ with  $\alpha=\{1,2,3\}$ forms a totally disconnected subposet matrix of $A.$ The submatrix  $D=A[{4}|{1,2,3}]$ has $3$ equal columns. By Theorem \ref{disconnectedsubposet}(i) the output poset matrix derived from $\square_i$ insertions at the labels $\{1,2,3\}$ are identical.}

\end{remark}

\begin{definition}
{\rm Let $\text{A}$ and $\text{B}$ be poset matrices of size $n\times n.$  If there exist a subposet matrix $\text{A}^{\prime\prime}$ of A which is a dual disconnected poset matrix of the  subposet matrix $\text{B}^{\prime\prime}$ of $\text{B}$ and the corresponding entries of $A$ and $B$ are always equal except at the region covered by the subposet matrices $\text{A}^{\prime\prime}$  and  $\text{B}^{\prime\prime}$ respectively,  then we shall henceforth refer to  the poset matrix $A$  as a \textbf{semi-equidual} of the poset matrix  $B$ and vice-versa.}
\end{definition}

\begin{example}{\rm Semi-equidual poset matrices of order $4$ and $5$ are as follows:}
\end{example}
$$A=\begin{blockarray}{ccccc}
1&2 & 3 &4\\
\begin{block}{(cccc)c}
 1 & 0 & 0&0 &1\\
 1 & 1 & 0&0 &2\\
 1 & 1 & 1&0 &3\\
 1 & 0 & 0&1 &4\\
\end{block}
\end{blockarray}\quad B=\begin{blockarray}{ccccc}
1&2 & 3 &4\\
\begin{block}{(cccc)c}
 1 & 0 & 0&0 &1\\
 1 & 1 & 0&0 &2\\
 1 & 0 & 1&0 &3\\
 1 & 0 & 1&1 &4\\
\end{block}
\end{blockarray}.$$

$$C=\begin{blockarray}{cccccc}
1&2 & 3 &4&5\\
\begin{block}{(ccccc)c}
 1 & 0 & 0&0&0 &1\\
 0 & 1 & 0&0& 0&2\\
 0 & 1 & 1&0 &0&3\\
 1&1&1&1&0&4\\
 1 & 1 & 1&1&1 &5\\
\end{block}
\end{blockarray}\quad D=\begin{blockarray}{cccccc}
1&2 & 3 &4&5\\
\begin{block}{(ccccc)c}
 1 & 0 & 0&0&0 &1\\
 1 & 1 & 0&0& 0&2\\
 0 & 0 & 1&0 &0&3\\
 1&1&1&1&0&4\\
 1 & 1 & 1&1&1 &5\\
\end{block}
\end{blockarray}.$$

\begin{remark}
{\rm The poset matrices $A$ and $B$ are semi-equidual poset matrices since as the disconnected subposet matrix $A[\{2,3,4\}]$  is a dual of the disconnected subposet matrix $B[\{2,3,4\}]$  and all entries of the column $1$  and row $1$ of $A$ are equal to the corresponding entries of $B$. Similarly, the poset matrices $C$ and $D$ are semi-equidual since as the disconnected subposet matrix $C[\{1,2,3\}]$ is a dual to the disconnected subposet matrix $D[\{1,2,3\}]$ and all the entries in the fourth row, fourth column, fifth row and fifth column of $C$ are equal to the corresponding entries of $D.$}
\end{remark}

\begin{theorem}\label{semiequidual}
 Let $A\in {\cal PM}(n)$ and let $\text{B}$ be a totally connected poset matrix of size $m\times  m.$Then the following holds.
\begin{enumerate}
\item If there exist a totally disconnected subposet matrix $A[\alpha]$ defined on $\alpha=\{1,...,k\}$ with $1<k<n$ such that the submatrix $D:=A[\{k+1,...,n\}|\{1,...,k\}]$ satisfies the condition that $D^{(p)}=D^{(p+1)}$ for each $p\in\{1,...,k-1\}.$ Then $\text{A}\square_1\text{B}$ \text{ is a semi-equidual poset matrix of} $\text{A}\square_k \text{B}$ and vice-versa.
\item If there exist a totally disconnected subposet matrix $A[\alpha]$ defined on $\alpha=\{k,...,n\}$ with $k\geq 2$ such that the submatrix $D:=A[\{k,...,n\}|\{1,...,k-1\}]$ satisfies the condition that $D_{(p)}=D_{(p+1)}$ for each $p\in\{1,...,n-k\}.$ Then $\text{A}\square_k\text{B}$ \text{ is a semi-equidual poset matrix of} $\text{A}\square_n \text{B}$ and vice-versa.
	\item If there exist a totally disconnected subposet matrix $A[\alpha]$ defined on $\alpha=\{d,...,k\}$ with $1<d<k<n$ such that the submatrix $D:=A[\{k,...,n\}|\{1,...,k-1\}]$ satisfies the condition that $D_{(p)}=D_{(p+1)}$ for each $p\in\{1,...,n-k\}$ where $D_{(p)}$ is a   \textbf{1}-vector or a \textbf{0}-vector of size $1 \times(k-1).$ Then $\text{A}\square_d\text{B}$ \text{ is a semi-equidual poset matrix of} $\text{A}\square_k \text{B}$ and vice-versa.
\end{enumerate}
\end{theorem}

\begin{proof}
Let $1,k\in\alpha$  as stated  in (1).  Consider that the output poset matrix structure of $A\square_{1}B=E$ can be partitioned into three submatrices $E1, E2$ and $E3$ such that:
\begin{itemize}
\item $E1=E[\{1,...,k+m-1\}]$  of size $(k+m-1)\times (k+m-1).$
\item $E2=E[\{k+m,...,n+m-1\}\mid\{1,...,k+m-1\}]$ of size $(n-k)\times (k+m-1).$
\item $E3=E[\{k+m,...,n+m-1\}]$ of size $(n-k)\times(n-k).$
\end{itemize}
\noindent $E1$ can be subpartitioned such that $E[\{1,...,m\}]=B$, $E[\{1+m,...,k+m-1\}\mid\{1,...,m\} ]$ is a zero matrix of size $(k-1)\times m,$ and $E[\{m+1,...,k+m-1\}]$ is either the connected poset matrix of size $1$ when $k=2$ or the totally disconnected poset matrix of size $k-1$ when $k>2.$

\noindent$E2$ can  be subpartitioned such that $E[\{k+m,...,n+m-1\}\mid\{1,...,m\}]={\mathbbm{1}}_{m}\otimes A^{(1)}[\{k+1,...,n\}\mid\{1,..,k\}]$,  and $E[\{k+m,...,n+m-1\}\mid\{m+1,...,k+m-1\}]=A[\{k+1,...,n\}\mid\{2,...,k\}].$

\noindent$E3=A[\{k+1,...,n\}].$

At the insertion point $k$, similarly consider that the output poset matrix structure of $A\square_{k}B=F$ can be partitioned into three submatrices $F1, F2$ and $F3$ such that:
\begin{itemize}
\item $F1=F[\{1,...,k+m-1\}]$  of size $(k+m-1)\times (k+m-1).$
\item $F2=F[\{k+m,...,n+m-1\}\mid\{1,...,k+m-1\}]$ of size $(n-k)\times (k+m-1).$
\item $F3=F[\{k+m,...,n+m-1\}]$ of size $(n-k)\times(n-k).$
\end{itemize}
\noindent $F1$ can be subpartitioned such that $F[\{1,...,k-1\}]$  is either the connected poset matrix of size $1$ when $k=2$ or the totally disconnected poset matrix of size $k-1$ when $k>2$
, $F[\{k,...,k+m-1\}\mid\{1,...,k-1\} ]$ is a zero matrix of size $m\times (k-1),$ and $F[\{k,...,k+m-1\}]=B.$

\noindent$F2$ can  be subpartitioned such that $F[\{k+m,...,n+m-1\}\mid\{k,...,k+m\}]={\mathbbm{1}}_{m}\otimes A^{(k)}[\{k+1,...,n\}\mid\{1,..,k\}]$, $F[\{k+m,...,n+m-1\}\mid\{1,...,k-1\}]=A[\{k+1,...,n\}\mid\{1,..,k-1\}].$

\noindent$F3=A[\{k+1,...,n\}].$

\noindent From the matrix structure of $F$ and $E,$ it follows that:
\begin{itemize}
\item For $E1=(e_{i,j})$ and $F1=(f_{i,j})$ both of size $k+m-1,$  it can be verified that $F1$ is a dual poset matrix of $E1$ since $F1=(e_{{k+m-j},{k+m-i}})$ and $E1=(f_{{k+m-j},{k+m-i}}).$ By definition \ref{disconnectMAT},  $F1$ and $E1$ are both disconnected poset matrix structure.
\item  $E2=F2$ since as each of the columns in the submatrix $A[\{k+1,...,n\}|\{1,...,k\}]$ has equal corresponding entries.
\item  $E3=F3.$
\end{itemize}

Thus, $E$ is a semi-equidual poset matrix of $F$ and vice-versa.

\noindent Let $k,n\in\alpha$  as stated  in (2).  Consider that the output poset matrix structure of $A\square_{k}B=G$ can be partitioned into three submatrices $G1, G2$ and $G3$ such that:
\begin{itemize}
\item $G1=G[\{1,...,k-1\}]$  of size $(k-1)\times (k-1).$
\item $G2=G[\{k,...,n+m-1\}\mid\{1,...,k-1\}]$ of size $(n+m-k)\times (k-1).$
\item $G3=G[\{k,...,n+m-1\}]$ of size $(n+m-k)\times(n+m-k).$
\end{itemize}
\noindent $G3$ can be subpartitioned such that $G[\{k,...,k+m\}]=B$, $G[\{k+m-1,...,n+m-1\}\mid\{k,...,k+m-1\} ]$  is a zero matrix of size $(n-k)\times m,$ and $G[\{k+m-1,...,n+m-1\}]$ is either the connected poset matrix of size $1$ when $k=2$ or the totally disconnected poset matrix of size $(n-k)\times(n-k)$ when $k>2.$

\noindent$G2$ can  be subpartitioned such that $G[\{k,...,k+m-1\}\mid\{1,...,k-1\}]={\mathbbm{1}}_m^T\otimes A_{(k)}[\{k,...,n\}\mid\{1,..,k-1\}]$,  and $G[\{k+m,...,n+m-1\}\mid\{1,...,k-1\}]=A[\{k,...,n\}\mid\{1,...,k-1\}].$

\noindent$G1=A[\{1,...,k-1\}].$

 At the insertion point $n,$ consider that the output poset matrix structure of $A\square_{n}B=H$ can be partitioned into three submatrices $H1, H2$ and $H3$ such that:
\begin{itemize}
\item $H1=H[\{1,...,k-1\}]$  of size $(k-1)\times (k-1).$
\item $H2=H[\{k,...,n+m-1\}\mid\{1,...,k-1\}]$ of size $(n+m-k)\times (k-1).$
\item $H3=H[\{k,...,n+m-1\}]$ of size $(n+m-k)\times(n+m-k).$
\end{itemize}
\noindent $H3$ can be subpartitioned such that $H[\{n,...,n+m-1\}]=B$, $H[\{n,...,n+m-1\}\mid\{k,...,n-1\} ]$ is the zero matrix of size $m\times (n-k),$ and $H[\{k,...,n-1\}]$ is either the connected poset matrix of size $1$ when $k=2$ or the totally disconnected poset matrix of size $(n-k)\times(n-k)$ when $k>2.$

\noindent$H2$ can  be subpartitioned such that $H[\{n,...,n+m-1\}\mid\{1,...,k-1\}]={\mathbbm{1}}_m^T\otimes A_{(n)}[\{k,...,n\}\mid\{1,..,k-1\}]$,  and $H[\{k,...,n-1\}\mid\{1,...,k-1\}]=A[\{k,...,n-1\}\mid\{1,...,k-1\}].$

\noindent$H1=A[\{1,...,k-1\}].$

\noindent From the matrix structure of $G$ and $H,$ it follows that:
\begin{itemize}
\item For $G3=(g_{i,j})$ and $H3=(h_{i,j})$ both of size $(n+m-k)\times(n+m-k),$  it can be verified that $H3$ is a dual poset matrix of $G3$ since $H3=(g_{{n+m-k+1-j},{n+m-k+1-i}})$ and $G3=(h_{{n+m-k+1-j},{n+m-k+1-i}}).$ By definition \ref{disconnectMAT},  $H3$ and $G3$ are disconnected poset matrix structures.
\item  $H2=G2.$ under the condition that each of the rows in the submatrix $A[\{k,...,n\}|\{1,...,k-1\}]$ has equal corresponding entries.
\item  $H3=G3.$
\end{itemize}
Thus,  $G$ is a semi-equidual poset matrix of $H$ and vice-versa.

The proof of (3) follows from the arguments in (1) and (2).

\end{proof}

\begin{example}
\end{example}
$\begin{blockarray}{ccccc}
1 & 2&3&4\\
\begin{block}{(cccc)c}
 1 & 0&0  &0&1\\
 1 & 1&0 &0&2\\
 1 & 0&1 &0&3 \\
 1&0&0&1&4\\
\end{block}
\end{blockarray}\scalebox{0.7}{$\square$}_2\;\;
\begin{blockarray}{ccc}
1 & 2 \\
\begin{block}{(cc)c}
 1 & 0&1 \\
 1 & 1&2\\
\end{block}
\end{blockarray}=
\begin{blockarray}{cccccc}
1&2 & 3 &4&5\\
\begin{block}{(ccccc)c}
 1 & 0 & 0&0&0 &1\\
 1 & 1 & 0&0& 0&2\\
 1 & 1 & 1&0 &0&3\\
 1&0&0&1&0&4\\
 1 & 0 & 0&0&1 &5\\
\end{block}
\end{blockarray}\quad\Leftrightarrow\quad$
\begin{minipage}[c]{.20\textwidth}
\begin{tikzpicture}
  [scale=.7,auto=center,every node/.style={circle,fill=blue!20}]

\node[rectangle] (c1) at (5.5,3) {2};
\node[rectangle] (c7) at (5.5,4.5) {3};
  \node[rectangle] (c2) at (7,1.5) {1};
  \node[rectangle] (c4) at (7,3) {4};
   \node[rectangle] (c5) at (8.5,3) {5};
      \draw (c5) -- (c2);
      \draw (c1) -- (c2);
        \draw (c2) -- (c4);
         \draw (c1) -- (c7);

 \end{tikzpicture}

\end{minipage}

$\begin{blockarray}{ccccc}
1 & 2&3&4\\
\begin{block}{(cccc)c}
 1 & 0&0  &0&1\\
 1 & 1&0 &0&2\\
 1 & 0&1 &0&3 \\
 1&0&0&1&4\\
\end{block}
\end{blockarray}\scalebox{0.7}{$\square$}_4\;\;
\begin{blockarray}{ccc}
1 & 2 \\
\begin{block}{(cc)c}
 1 & 0&1 \\
 1 & 1&2\\
\end{block}
\end{blockarray}=
\begin{blockarray}{cccccc}
1&2 & 3 &4&5\\
\begin{block}{(ccccc)c}
 1 & 0 & 0&0&0 &1\\
 1 & 1 & 0&0& 0&2\\
 1 & 0 & 1&0 &0&3\\
 1&0&0&1&0&4\\
 1 & 0 & 0&1&1 &5\\
\end{block}
\end{blockarray}\quad\Leftrightarrow\quad$
\begin{minipage}[c]{.20\textwidth}
\begin{tikzpicture}
  [scale=.7,auto=center,every node/.style={circle,fill=blue!20}]

\node[rectangle] (c1) at (5.5,3) {2};
\node[rectangle] (c7) at (8.5,4.5) {5};
  \node[rectangle] (c2) at (7,1.5) {1};
  \node[rectangle] (c4) at (7,3) {3};
   \node[rectangle] (c5) at (8.5,3) {4};
      \draw (c5) -- (c2);
      \draw (c1) -- (c2);
        \draw (c2) -- (c4);
         \draw (c5) -- (c7);

 \end{tikzpicture}

\end{minipage}

 \begin{remark}
{\rm By applying Theorem  \ref{semiequidual}(ii), we can observe that  the input disconnected subposet matrix $A[ \{2,3,4 \}]$ forms semi-equidual poset matrices at insertion points $2$ and $4,$ and we can observe that the submatrix  $D=A[ \{2,3,4 \}| \{1 \}]$ has all equal rows.}
 \end{remark}

\begin{example}
\end{example}
$\begin{blockarray}{ccccc}
1 & 2&3&4\\
\begin{block}{(cccc)c}
 1 & 0&0  &0&1\\
 0 & 1&0 &0&2\\
 1 & 1&1 &0&3 \\
 1&1&1&1&4\\
\end{block}
\end{blockarray}\scalebox{0.7}{$\square$}_1\;\;
\begin{blockarray}{ccc}
1 & 2 \\
\begin{block}{(cc)c}
 1 & 0&1 \\
 1 & 1&2\\
\end{block}
\end{blockarray}=
\begin{blockarray}{cccccc}
1&2 & 3 &4&5\\
\begin{block}{(ccccc)c}
 1 & 0 & 0&0&0 &1\\
 1 & 1 & 0&0& 0&2\\
 0 & 0 & 1&0 &0&3\\
 1&1&1&1&0&4\\
 1 & 1 & 1&1&1 &5\\
\end{block}
\end{blockarray}\quad\Leftrightarrow\quad$
\begin{minipage}[c]{.20\textwidth}
\begin{tikzpicture}
  [scale=.7,auto=center,every node/.style={circle,fill=blue!20}]

\node(d4) at (10,2) {5};
\node(d1) at (10,0.5) {4};
\node(d6) at (11,-0.3)  {2};
\node(d9) at (11,-1.5)  {1};
\node(d7) at (9,-0.3) {3};
\draw (d6) -- (d1);
\draw (d1) -- (d7);
\draw (d1) -- (d4);
\draw (d9) -- (d6);

 \end{tikzpicture}

\end{minipage}

$\begin{blockarray}{ccccc}
1 & 2&3&4\\
\begin{block}{(cccc)c}
 1 & 0&0  &0&1\\
 0 & 1&0 &0&2\\
 1 & 1&1 &0&3 \\
 1&1&1&1&4\\
\end{block}
\end{blockarray}\scalebox{0.7}{$\square$}_2\;\;
\begin{blockarray}{ccc}
1 & 2 \\
\begin{block}{(cc)c}
 1 & 0&1 \\
 1 & 1&2\\
\end{block}
\end{blockarray}=
\begin{blockarray}{cccccc}
1&2 & 3 &4&5\\
\begin{block}{(ccccc)c}
 1 & 0 & 0&0&0 &1\\
 0 & 1 & 0&0& 0&2\\
 0 & 1 & 1&0 &0&3\\
 1&1&1&1&0&4\\
 1 & 1 & 1&1&1 &5\\
\end{block}
\end{blockarray}\quad\Leftrightarrow\quad$
\begin{minipage}[c]{.20\textwidth}
\begin{tikzpicture}
  [scale=.7,auto=center,every node/.style={circle,fill=blue!20}]

\node(d4) at (10,2) {5};
\node(d1) at (10,0.5) {4};
\node(d6) at (11,-0.3)  {1};
\node(d7) at (9,-0.3) {3};
\node(d9) at (9,-1.5) {2};
\draw (d6) -- (d1);
\draw (d1) -- (d7);
\draw (d1) -- (d4);
\draw (d9) -- (d7);

 \end{tikzpicture}

\end{minipage}

 \begin{remark}
{\rm By applying Theorem  \ref{semiequidual}(i), we can observe that  the input disconnected subposet matrix $A[ \{1,2 \}]$ forms semi-equidual poset matrices at insertion points $1$ and $2,$ and we can observe that the submatrix  $D=A[ \{3,4 \}| \{1,2 \}]$ is comprised of equal columns.}
 \end{remark}

\end{document}